%% file: paper.tex
\documentclass[reqno]{siamart}

\usepackage{graphicx}
\usepackage[T1]{fontenc}
\usepackage{amsmath}
\usepackage{amssymb}
\usepackage{stmaryrd}

\newcommand{\dd}{\, \mathrm{d}}
\newcommand{\E}{\mathcal{E}}

\newcommand{\M}{\mathcal{M}}
\newcommand{\md}{\partial^\bullet}
\newcommand{\meas}{\mathrm{meas}}
\newcommand{\R}{\mathbb{R}}
\newcommand{\T}{\mathcal{T}}

\newcommand{\bigjump}[1]{\bigl\llbracket #1 \bigr\rrbracket}
\newcommand{\bigaverage}[1]{\bigl\{ #1 \bigr\}}
\newcommand{\old}{\text{old}}

\newtheorem{remark}{Remark}
\newtheorem{scheme}{Scheme}

\crefformat{subsection}{#2Section~#1#3}
\Crefformat{subsection}{#2Section~#1#3}

\title{An unfitted discontinuous Galerkin scheme for conservation laws on evolving surfaces%
  \thanks{The work was partially supported by the German Research
    Foundation (DFG), grant EN 1042/4-1 and by the UK Engineering and
    Physical Sciences Research Council, grant EPSRC EP/J004057/1.}  }

\author{Christian Engwer\thanks{Institute for Computational and Applied Mathematics, University of M\"{u}nster, Germany.
    (\email{christian.engwer@uni-muenster.de},
    \email{sebastian.westerheide@uni-muenster.de})} %
  \and Thomas Ranner\thanks{School of Computing, EC Stoner Building, University of Leeds, UK. (\email{T.Ranner@leeds.ac.uk})} %
  \and Sebastian Westerheide\footnotemark[2]}

\headers{UDG for conservation laws on evolving surfaces}{C.\ Engwer, T.\ Ranner and S.\ Westerheide}

\ifpdf
\hypersetup{
  pdftitle={UDG for conservation laws on evolving surfaces},
  pdfauthor={C.\ Engwer, T.\ Ranner and S.\ Westerheide}
}
\fi

\begin{document}

\mbox{}

\maketitle

\begin{abstract}
  Motivated by considering partial differential equations arising from conservation laws posed on evolving surfaces, a new numerical method for an advection problem is developed and simple numerical tests are performed.
  The method is based on an unfitted discontinuous Galerkin approach where the surface is not explicitly tracked by the mesh which means the method is extremely flexible with respect to geometry.
  Furthermore, the discontinuous Galerkin approach is well-suited to capture the advection driven by the evolution of the surface without the need for a space-time formulation, back-tracking trajectories or streamline diffusion.
  The method is illustrated by a one-dimensional example and numerical results are presented that show good convergence properties for a simple test problem.
\end{abstract}

\begin{keywords}
  Discontinuous Garlerkin; unfitted finite elements; surface partial differential equations; conservation laws.
\end{keywords}

\begin{AMS}
35L02, 
35L65, 
35Q90, 
35R01, 
65N12, 
65N30  
\end{AMS}

\section{Introduction}

Interest in the study of partial differential equations on evolving surfaces has grown in recent years with applications in many areas including materials science (e.g.\ diffusion of species along grain boundaries \cite{DecEllSty01}), fluid dynamics (e.g.\ surface active agents along the interface between two fluids \cite{JamLow04}) and cell biology (e.g.\ cell motility involving the processes on the cell membrane \cite{NeiMacWeb11}), for example. See the review of \cite{DziEll13a} for a more detailed list.
In this work, we derive a new numerical scheme for an essential conservation law on evolving hypersurfaces and present first numerical results.

\subsection{Problem}

As a model problem, we consider the evolution of a conserved quantity on an evolving curve or surface.
Fix $T > 0$ and let $\{ \Gamma(t) \}_{t \in [0,T]}$ be a time-dependent,
connected, compact,
smooth $n$-dimensional hypersurface embedded in $\R^{n+1}$ for
$n \in \{ 1,2 \}$, $\vec{\nu}( \cdot, t ) \colon \Gamma(t) \to \R^{n+1}$ denote
a field of unit normal vectors to $\Gamma(t)$, and write $\Gamma_0 = \Gamma(0)$.
We consider two different descriptions of the surface.
First, we say $\Gamma(t)$ is defined by a diffeomorphic parametrization $G( \cdot, t ) \colon \Gamma_0 \to \R^{n+1}$ so that $\Gamma(t) = G( \Gamma_0, t )$ for $t \in [0,T]$.
The map $G$ defines a field $\vec{w}$ which describes the velocity of $\Gamma(t)$,
precisely by $\partial_t G( \cdot, t ) = \vec{w}( G( \cdot, t ), t )$.
Second, we describe $\Gamma(t)$ as the zero level set of a smooth function $\Phi( \cdot, t ) \colon \R^{n+1} \to \R$ so that $\Gamma(t) = \{ x \in \R^{n+1} : \Phi( x, t ) = 0 \}$.
We will use the diffeomorphic parametrization to define the equations we solve and the level set function to define the geometry in our computational scheme.

Consider a control volume $\M(t) \subset \Gamma(t)$ which is the image of a control volume $\M_0 \subset \Gamma_0$ under the flow $G$, i.e.\ $\M(t) = G( \M_0, t )$.
Let $u$ denote a time-dependent field on $\Gamma(t)$ which satisfies a conservation law of the form
\begin{equation}
  \label{eq:balance}
  \frac{\mathrm{d}}{\mathrm{d}t} \int_{\M(t)} u \dd \sigma
  = - \int_{\partial \M(t)} \vec{q} \cdot \vec{\mu} \dd \xi,
\end{equation}
where $\vec{q}$ is a tangential flux on $\Gamma(t)$ and $\vec{\mu}$ is the co-normal vector to $\partial \M(t)$ ($\vec{\mu}$ is normal to $\partial \M(t)$ and $\vec{\mu} \cdot \vec{\nu} = 0$).
Applying a transport formula to the left hand side of \cref{eq:balance} and the divergence theorem for hypersurfaces to the right hand side yields
\begin{equation}
  \label{eq:local-cons}
  \int_{\M(t)} \md u + u \nabla_\Gamma \cdot \vec{w} + \nabla_\Gamma \cdot \vec{q} \dd \sigma = 0,
\end{equation}
where $\md u$ is the material derivative of $u$ and $\nabla_\Gamma \cdot {}$ denotes the tangential divergence operator. See \cref{sec:prelim} for details.
We describe equation \cref{eq:local-cons} as a \emph{local} conservation law for $u$ or a \emph{global} conservation law in the case $\M(t) = \Gamma(t)$.

Since \cref{eq:local-cons} holds for any
$\M(t)$, we
arrive at a pointwise conservation law.
We wish to find a time-dependent surface field
$u \colon \bigcup_{t \in [0,T]} \Gamma(t) \times \{ t \} \to \R$ with
\begin{equation}
  \label{eq:problem}
  \md u + u \nabla_\Gamma \cdot \vec{w} + \nabla_\Gamma \cdot \vec{q} = 0 \qquad \mbox{ on } \bigcup_{t \in [0,T]} \Gamma(t) \times \{ t \}.
\end{equation}
For more details on the notation, see \cite{DziEll13a}.

In this work, we aim to derive a new numerical method using an implicit representation of the surface for which we can show a discrete analogue to the global conservation law.
To solve problems of this form, we are required to consider two separate problems. We have to discretize the surface partial differential equation at each given time step and we have to evaluate the transport from the old surface to the new surface.
Either these two components are computed in a splitting approach and discretized separately,
or computed directly in a coupled manner.
In either case, we must understand the effects of the flux $\vec{q}$ and the advection of the moving surface separately.
Previously (see the following section), attention has been paid to the stationary case and the treatment of the flux $\vec{q}$.
In this work, we are interested in the advection driven by the moving surface.
We thus restrict to the case $\vec{q} = 0$ and a single time step, which simplifies considerations.
Extensions to the truly time-dependent case and $\vec{q} \neq 0$ will be considered in future work.

\subsection{Previous approaches}

Previous computational approaches lie in two categories. The first is based on a moving mesh approach and the second is based on a static mesh.

For surface PDEs, the basic ideas of the moving mesh approach date back to the ideas of the
finite element method for elliptic PDEs on stationary surfaces which was
presented in \cite{Dzi88}.
A smooth surface is approximated by a union of elements (usually simplices but also possibly quadradrilaterals) whose vertices
move according to a globally defined smooth velocity.
It was applied by \cite{DziEll07}, using a finite element discretization
for an advection-diffusion problem, and by \cite{DziKroMul13}, using a
finite volume method for conservation laws.
Moving mesh schemes have proven to be useful in many practical situations
and can be shown analytically to be stable and accurate, or preserve various conservative properties.
Their main disadvantage is that the moving meshes may degenerate, leading
to frequent remeshing \cite{EilEll08}.
Extensions have been presented recently by \cite{EllSty12,EllVen15} that
reduce the need for remeshing by allowing mesh points to move with a different
non-physical velocity.
However, in many situations preserving a good mesh remains very challenging.

In the static mesh approach, an arbitrary fixed background mesh is used and the evolution of the surface is defined implicitly.
One approach is to use a level set function together with an unfitted finite element method where the computational domain consists of partial cut-cell elements \cite{DecDziEll09,OlsReuGra09,DecEllRan14,OlsReuXu14a}.
The authors of \cite{OlsReuXu14} apply a space-time formulation to produce a stable, accurate method, but the method is only approximately globally conservative and requires the use of complicated four-dimensional space-time elements which are not available for geometries other than simplices.
An alternative method proposed by \cite{DecEllRan14} uses a semi-Lagrangian formulation that employs a non-local right hand side.
The resulting scheme recovers a globally mass conservation property and performs well in practice for a test problem on a curve, but lacks analysis and the non-local term is expensive to compute, especially in three space dimensions.
A phase field representation of the interface was used by \cite{EllStiSty11,TeiLiLow09}.
The authors of \cite{EllStiSty11} use a narrow band formulation and recover a discrete analogue to the continuous level global conservation law but must use stream line diffusion
to stabilize the scheme.

A key difficulty to overcome with the static mesh approach is to derive a stable
scheme that adequately treats the advection driven by the evolution of the
surface.
The moving mesh approach implicitly deals with this problem by using moving basis functions.
The advective flux usually has a component orthogonal to the surface.
Therefore, any other fluxes (e.g.\ diffusive surface fluxes) can not
stabilize the method.
Furthermore, the previous ideas based on finite element methods are not
well-suited to advection-dominated problems.
Extra complications such as space-time formulations, semi-Lagrangian
terms or streamline diffusion are required.
In this work, we use a more suitable discontinuous Galerkin (DG) method
which can naturally handle advection or advection-dominated problems.

\section{Computational approach}

Our approach is based on reformulating problem \cref{eq:problem} as a sequence
of bulk advection problems with singular source and sink terms and applying the
unfitted discontinuous Galerkin (UDG) method \cite{BasEng09} on a static bulk mesh.
This mesh can be chosen independent of the evolving surface, the surface is not
explicitly represented by the mesh.
Our formulation guarantees a globally conservative scheme.

\subsection{Preliminaries}
\label{sec:prelim}

Let ${U}$ be a polygonal domain in $\R^{n+1}$ which contains $\Gamma(t)$ for all times $t \in [0,T]$.
We suppose that for each time $t \in [0,T]$, $\Gamma(t)$ is defined as the zero level set of a smooth function $\Phi( \cdot, t ) \colon {U} \to \R$ and assume $\nabla \Phi( x, t ) \neq 0$ for $x \in {U}$, $t \in [0,T]$.
Let $\vec{\nu}$ be the field of unit normal vectors to $\Gamma(t)$ oriented by $\vec{\nu} = \nabla \Phi / | \nabla \Phi |$.
The normal component of the velocity field is given by
$\vec{w} \cdot \vec{\nu} = -{ \Phi_t }/{ | \nabla \Phi | }$.

The tangential gradient of a smooth surface field $\eta \colon \Gamma(t) \to \R$ can be defined as
\begin{equation*}
  \nabla_\Gamma \eta := \nabla \tilde\eta - \bigl( \nabla \tilde\eta \cdot \vec{\nu} \bigr) \vec{\nu},
\end{equation*}
where $\tilde\eta$ is a spatially differentiable extension of $\eta$ to a neighborhood of $\Gamma(t)$
and $\nabla \tilde\eta$ is its gradient with respect to the ambient Cartesian coordinates.
The tangential gradient $\nabla_\Gamma \eta$ has $n+1$ components $( \underbar{D}_1 \eta, \ldots, \underbar{D}_{n+1} \eta )$.
For a smooth surface vector field $\vec{v}$, the tangential divergence is given by
\begin{equation*}
  \nabla_\Gamma \cdot \vec{v} = \sum_{i=1}^{n+1} \underbar{D}_i \vec{v}_i.
\end{equation*}

The material derivative of a smooth, time-dependent surface field $\eta( \cdot, t ) \colon \Gamma(t) \to \R$ is defined as
\begin{equation*}
  \md \eta := \partial_t \tilde\eta + \vec{w} \cdot \nabla \tilde\eta,
\end{equation*}
where again $\partial_t \tilde\eta$ and $\nabla \tilde\eta$ are Cartesian derivatives of a differentiable extension $\tilde\eta$ of $\eta$ to a space-time neighborhood of $\{ \Gamma(t) \}_{t \in [0,T]}$.
The material derivative describes the variation of $\eta$ with respect to the evolution of the surface.
It also has a key role in the following transport relation for evolving material surfaces, see e.g.\
\cite{CerFriGur05,DziEll07}:
\begin{equation}
  \label{eq:transport}
  \frac{\mathrm{d}}{\mathrm{d}t} \int_{\Gamma(t)} \eta \dd \sigma = \int_{\Gamma(t)} \md \eta + \eta \nabla_\Gamma \cdot \vec{w} \dd \sigma.
\end{equation}

\subsection{Motivation}
\label{sec:motivation}

Proceeding formally,
we multiply problem \cref{eq:problem} (with $\vec{q}=0$) by a smooth function $\varphi \colon U \to \R$, integrate over $\Gamma(t)$ for $t \in [0, T]$ and apply transport relation \cref{eq:transport} to see
\begin{align*}
  \frac{\mathrm{d}}{\mathrm{d}t} \int_{\Gamma(t)} u \varphi \dd \sigma
  - \int_{\Gamma(t)} u \md \varphi \dd \sigma = 0.
\end{align*}
Since $\varphi$ does not depend on time, we have $\md \varphi = \vec{w} \cdot \nabla \varphi$.

\begin{figure}[tb]
  \centering
  \includegraphics[width=0.44\linewidth]{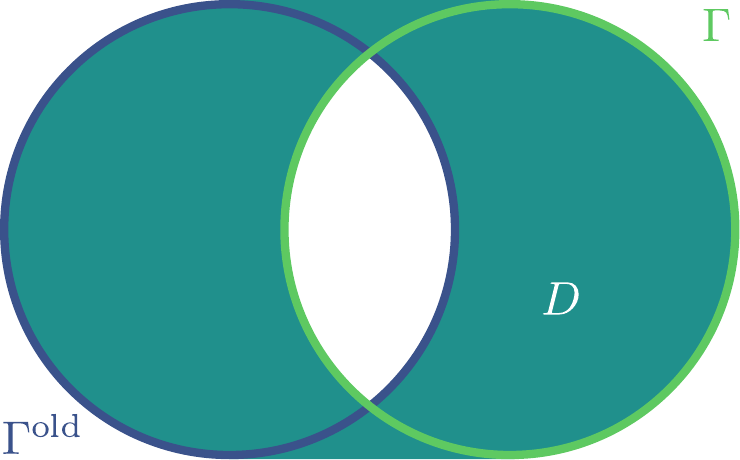}
  \caption{An example of the geometry configuration.
           Here, a circle is translated with a horizontal velocity between time steps.
           Note that a time step with a large time step size $\tau$ is illustrated.}
  \label{fig:geometry}
\end{figure}

We fix a time $t^* \in [0, T]$ and $\tau > 0$ small enough that $t^* - \tau \ge 0$.
Integrating in time over the time interval $[ t^* - \tau, t^* ]$ yields
\begin{equation*}
  \int_{\Gamma(t^*)} u(t^*) \varphi \dd \sigma - \int_{\Gamma(t^* - \tau)} u(t^* - \tau) \varphi \dd \sigma
  - \int_{t^* - \tau}^{t^*} \left( \int_{\Gamma(t)} u \vec{w} \cdot \nabla \varphi \dd \sigma \right) \dd t= 0.
\end{equation*}
We will approximate the space-time integral in the third term by a scaled bulk integral over the spatial-only domain
\begin{equation*}
  D := \bigcup\nolimits_{t \in [ t^* - \tau, t^* ]} \Gamma(t) \, ,
\end{equation*}
which corresponds to the projection of the space-time domain
to spatial-only coordinates.
An example of the geometric setup is shown in \cref{fig:geometry}.

Denoting an extension of a function $\psi(\cdot,t^*)$ by $\psi^e \colon D \to \R$,
we have
\begin{alignat*}{1}
  - \int_{t^* - \tau}^{t^*} &\left( \int_{\Gamma(t)} u \vec{w} \cdot
                              \nabla \varphi \dd \sigma \right) \dd t
  \approx -\tau \int_{\Gamma(t^*)} u \vec{w} \cdot \nabla \varphi \dd \sigma \\
  & \quad \approx -\frac{\tau}{\gamma} \int_{-\gamma^-}^{\gamma^+}
    \left(
    \int_{\Gamma_s(t^*)} u^e \vec{w}^e \cdot \nabla
    \varphi \dd \sigma \right)
    \dd s
  = - \frac{\tau}{\gamma}
    \int_{D} u^e \vec{w}^e \bigl| \nabla \Phi(\cdot,t^*) \bigr| \cdot \nabla \varphi \dd x,
\end{alignat*}
where we employ the right-hand rectangle method,
an approximation of the resulting surface integral using the level sets
$\Gamma_s(t^*) := \{ x \in D : \Phi(x,t^*) = s \}$
of $\Phi(\cdot,t^*) \big\rvert_{D}$ and apply the coarea formula,
while assuming that the time step size $\tau$ as well as
$\gamma := \gamma^- + \gamma^+$ with $\gamma^- := - \min_D{\Phi(\cdot,t^*)}$
and $\gamma^+ := \max_D{\Phi(\cdot,t^*)}$ are small.
Note that the assumption on $\gamma$ is reasonable for a small time step size.
If $\vec{w}$ is constant and $\Phi$ is a signed distance function, i.e.\
$| \nabla \Phi | \equiv 1$, $\gamma$ corresponds to the travel distance of an
individual point on the surface.

In order to be able to discretize using flux-based numerical schemes, such as
unfitted DG schemes, we use integration by parts to reformulate this integral
via the divergence of an advective flux tested with $\varphi$.
As an approximate reformulation of problem \cref{eq:problem} over the time
interval $[t^* - \tau, t^*]$, we therefore obtain the following stationary
problem:
seek $u^e \colon D \to \R$ with
\begin{multline}
  \label{eq:problem_reformulated}
  \int_{\Gamma} u^e \varphi \dd \sigma - \int_{\Gamma^{\old}} u^{\old} \varphi \dd \sigma \\
  + \frac{\tau}{\gamma} \int_{D} \nabla \cdot \Big( u^e \vec{w}^e \bigl| \nabla \Phi(\cdot,t^*) \bigr| \Big) \varphi \dd x
  - \frac{\tau}{\gamma} \int_{\partial D} u^e \vec{w}^e \bigl| \nabla \Phi(\cdot,t^*) \bigr| \cdot \nu_D \varphi \dd \sigma
  = 0,
\end{multline}
for all smooth functions $\varphi \colon D \to \R$.
Here and in the following, we use the notation $\Gamma = \Gamma(t^*)$,
$\Gamma^{\old} = \Gamma(t^* - \tau)$ and $u^{\old} = u(t^* - \tau)$,
and $\nu_D$ denotes the outward pointing unit normal vector field to $\partial D$.

\smallskip

\begin{remark}
  In the special case
  $\partial D = \Gamma^{\old} \mathbin{\dot{\cup}} \Gamma$
  (e.g.\ a geometrical setup as in \cref{fig:domain-ex}),
  problem \cref{eq:problem_reformulated} is consistent with
  the following bulk advection problem in strong form:
  \vspace*{-0.5\baselineskip}
  \begin{alignat*}{3}
    \nabla \cdot \Big( u^e \vec{w}^e \bigl| \nabla \Phi(\cdot,t^*) \bigr| \Big) &= 0 \quad & &\text{in} \ \ & &D, \\
    u^e \vec{w}^e \bigl| \nabla \Phi(\cdot,t^*) \bigr| \cdot \nu_D &= - \frac{\gamma}{\tau} u^{\old} \quad & &\text{on} \ \ & &\Gamma^{\old}, \\
    u^e \vec{w}^e \bigl| \nabla \Phi(\cdot,t^*) \bigr| \cdot \nu_D &= \frac{\gamma}{\tau} u^e \quad & &\text{on} \ \ & &\Gamma. \\
  \end{alignat*}
\end{remark}

\vspace*{-1.5\baselineskip}

\subsection{Unfitted discrete geometry}
\label{sec:discrete_geometry}

\begin{figure}[tb]
  \centering
  \includegraphics[width=0.325\textwidth]{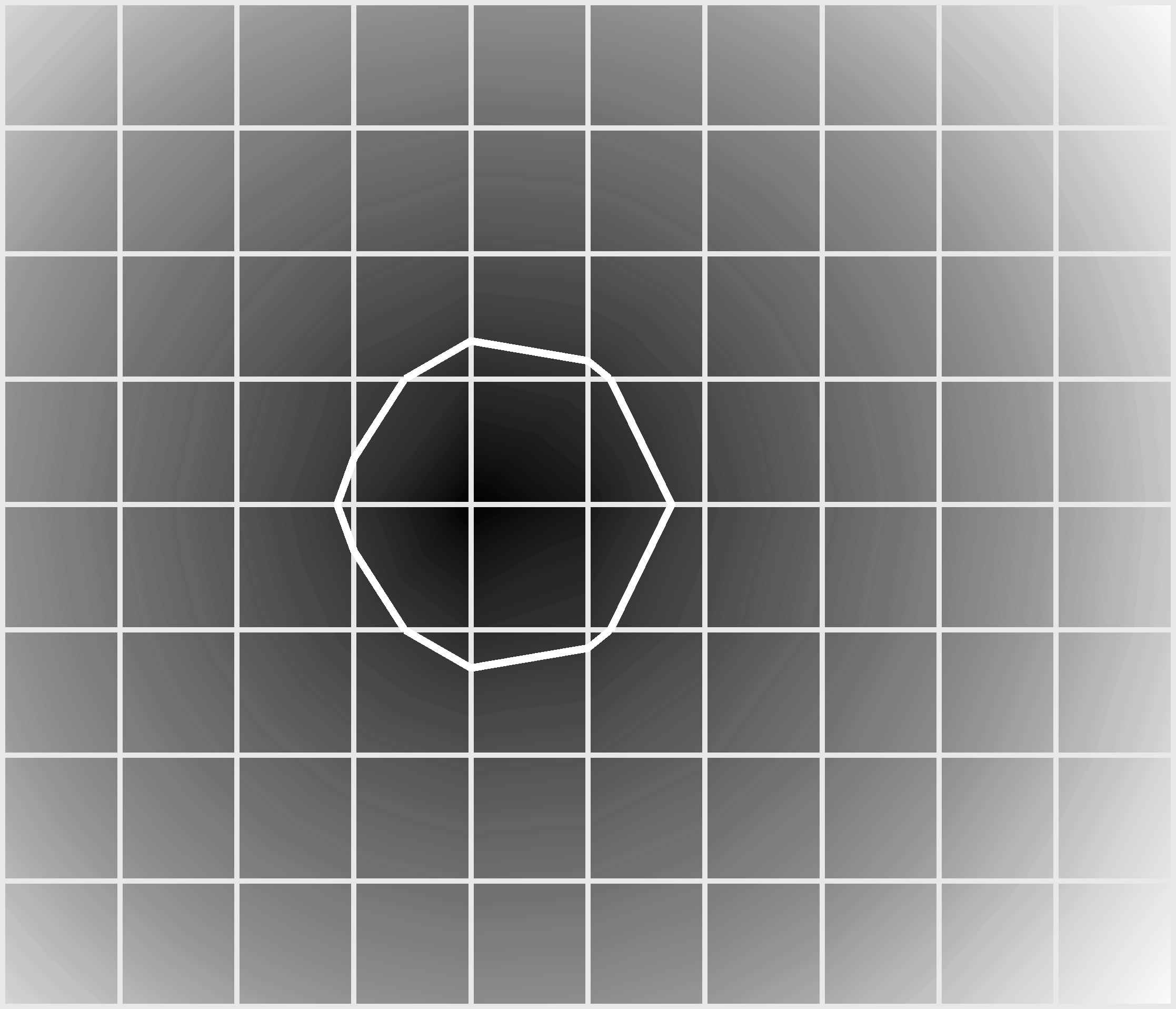}\qquad
  \includegraphics[width=0.325\textwidth]{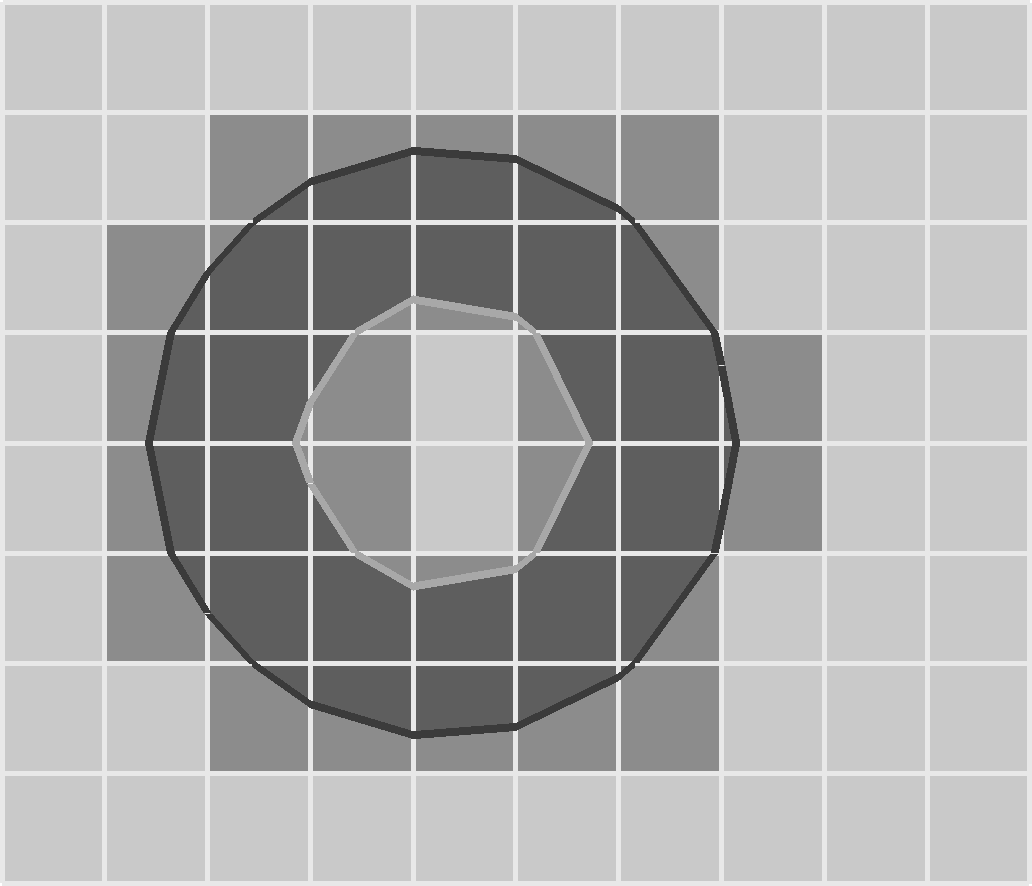}
  \caption{Example of an unfitted discrete geometry.
    The left image shows a discrete level set function $\Phi_h$ and its zero level set in white when $\Gamma$ is a circle.
    The right image shows an example of $\tilde{\T}_h$, $\Gamma_h$, $\T_h$, the cut cell mesh $\hat{\T}_h$ for computational domain $D_h$, and $\Gamma_h^{\old}$ (from light gray to dark gray, where pixels which represent meshes have the color of the most specialized mesh) when the initial curve $\Gamma^{\old}$ is shrunk to a smaller circle $\Gamma$.
  }
  \label{fig:domain-ex}
\end{figure}

In order to discretize problem \cref{eq:problem_reformulated}, we start by discretizing the geometry.
Let $\tilde{\T}_h$ be a shape regular decomposition of ${U}$ into closed elements, either tetrahedra or hexahedra for $n=2$, triangles or quadrilaterals for $n=1$, and denote by $h$ the maximum element size.
Let $X_h$ be the space of piecewise linear (for simplices), bilinear (for quadrilaterals) or trilinear (for hexahedra) continuous functions over $\tilde{\T}_h$ and denote by $I_h$ interpolation of functions in $C({U})$ into $X_h$.
We will write $\Phi_h( \cdot, t ) = I_h \Phi( \cdot, t )$ for a discrete level set function and set
\begin{align*}
  \Gamma_h & := \{ x \in {U} : \Phi_h( x, t^* ) = 0 \}, \quad
  \Gamma_h^{\old} := \{ x \in {U} : \Phi_h( x, t^* - \tau ) = 0 \} \\
  D_h & := \{ x \in {U} : \text{there exists } t \in [t^* - \tau, t^*] \mbox{ with } \Phi_h( x, t ) = 0 \}.
\end{align*}
Note that $\Gamma_h^{\old}, \Gamma_h \subset D_h$ since $D_h$ is a closed set.
We also use the notation
\begin{align*}
  \T_h := \left\{ K \in \tilde\T_h : \meas( K \cap D_h ) > 0 \right\}.
\end{align*}
An example is shown in \cref{fig:domain-ex}.
It will be useful to consider the restriction of the decomposition $\T_h$ to the computational domain $D_h$:
\begin{align*}
  \hat{\T}_h := \left\{ \hat{K} = K \cap D_h : K \in \T_h \right\}.
\end{align*}
We note that the $\hat{K} \in \hat{\T}_h$ are arbitrary shaped elements and
call those elements cut cells.
In general, they are not either shape regular or even convex.
We will also consider the internal faces of the cut cell mesh $\hat{\T}_h$,
denoted by $\hat{\E}_h$.
The set $\hat{\E}_h$ is often called the skeleton of the mesh.
To each internal face $\hat{E} \in \hat{\E}_h$, which is the intersection of two elements $\hat{K}^+, \hat{K}^- \in \hat{\T}_h$, we assign a unit normal vector field $\vec{\nu}_{\hat{K}^+} = - \vec{\nu}_{\hat{K}^-}$ and arbitrarily choose $\vec{\nu}_{\hat{E}} = \vec{\nu}_{\hat{K}^+}$.

\medskip

\begin{remark}
  Our solution will be defined over the union of elements in $\T_h$,
  but we are only interested in the values of our solution variable over the sharp interface approximation $\Gamma_h$.
  The integrals in the method will be computed over either the sharp interface approximation $\Gamma_h$ or the unfitted cut cell mesh $\hat{\T}_h$.
  This is a similar approach to \cite{DecEllRan14}, but different to \cite{OlsReuGra09}, thus avoiding difficulties in defining our discrete spaces and constructing a natural basis.
\end{remark}

\smallskip

\subsection{Method}
\label{sec:method}

We introduce the discrete spaces $V_h$ given by
\begin{align*}
  V_h := \bigl\{ \varphi_h = ( \varphi_K )_{K \in \T_h} : \varphi_K \in P^k( K ) \mbox{ for all } K \in \T_h \bigr\},
\end{align*}
where $P^k( K )$ denotes the space of piecewise polynomials of degree $k$ over the volumetric domain $K$.
Note that functions in $V_h$ are discontinuous and do not take a unique value along the faces $\hat{\E}_h$, in general.
On each face $\hat{E} \in \hat{\E}_h$ with adjacent elements $\hat{K}^+, \hat{K}^- \in \hat{\T}_h$ as defined above,
we define the jump of a function $\varphi_h \in V_h$ as $\bigjump{ \varphi_h } := \varphi_h \rvert_{\hat{K}^+} - \varphi_h \rvert_{\hat{K}^-}$
and its average as $\bigaverage{ \varphi_h } := \frac{1}{2} \bigl( \varphi_h \rvert_{\hat{K}^+} + \varphi_h \rvert_{\hat{K}^-} \bigr)$.

We discretize problem \cref{eq:problem_reformulated} by approximating
$u^e$ by $u_h \in V_h$, the test functions $\varphi$ by $\varphi_h \in V_h$,
$\Phi$ by $\Phi_h$ and $\vec{w}^e$ by a discrete approximation $\vec{w}_h$
(which is given in more detail in the sequel).
Integrals over $\Gamma$ and $\Gamma^{\old}$ are approximated by integrals
over $\Gamma_h$ and $\Gamma_h^{\old}$, respectively, and integrals over $D$
are approximated by integrals over $D_h$.
Finally, we split the bulk integral into a sum of integrals over the cut
cells $\hat{K} \in \hat{\T}_h$ and integrate by parts on each $\hat{K}$.
Using classical upwind stabilization (as in \cite{engwer2010:transport}),
we obtain the following unfitted DG scheme:

\medskip

\begin{scheme}
  Given $u_h^{\old} \in L^2(\Gamma_h^{\old})$, find $u_h \in V_h$, such that
  \begin{multline}
    \label{eq:scheme}
    \int_{\Gamma_h} u_h \varphi_h \dd \sigma_h
    + \frac{\tau}{\gamma} \bigg( \sum_{\hat{E} \in \hat{\E}_h} \int_{\hat{E}} u_h^{\uparrow} \bigjump{\varphi_h} \bigaverage{\vec{w}_h | \nabla \Phi_h(\cdot,t^*) |} \cdot \vec{\nu}_{\hat{E}} \dd \sigma_h \\
      - \sum_{\hat{K} \in \hat{\T}_h} \int_{\hat{K}} u_h \vec{w}_h | \nabla \Phi_h(\cdot,t^*) | \cdot \nabla \varphi_h \dd x_h \bigg) =
    \int_{\Gamma_h^{\old}} u_h^{\old} \varphi_h \dd \sigma_h
    \quad \text{ for all } \varphi_h \in V_h,
  \end{multline}
  where $u_h^{\uparrow}$ denotes the upwind solution of $u_h$ on $\hat{\E}_h$ which is given by
  \begin{align*}
    u_h^{\uparrow} \bigr\rvert_{\hat{E}} =
    \begin{cases}
      u_h \rvert_{\hat{K}^+} & \mbox{ if } \bigaverage{\vec{w}_h | \nabla \Phi_h(\cdot,t^*) |} \cdot \vec{\nu}_{\hat{E}} \ge 0 \\
      u_h \rvert_{\hat{K}^-} & \mbox{ if } \bigaverage{\vec{w}_h | \nabla \Phi_h(\cdot,t^*) |} \cdot \vec{\nu}_{\hat{E}} < 0
    \end{cases}
    \qquad
    \text{ for each } \hat{E} \in \hat{\E}_h.
  \end{align*}

  Note that if a section of $\Gamma_h$ or $\Gamma_h^{\old}$ is part of
  a face $\hat{E} \in \hat{\E}_h$, we choose $u_h$ as $u_h^\uparrow$
  and $u_h^{\old}$ as $u_h^{\old,\uparrow}$, respectively.
\end{scheme}

\vspace{0.10\baselineskip}

\begin{remark}
  We note that, since this construction implies that $\varphi_h \equiv 1$ is
  an admissible test function, a global mass conservation law is recovered:
  \begin{align*}
    \int_{\Gamma_h} u_h \dd \sigma_h
    = \int_{\Gamma_h^{\old}} u_h^{\old} \dd \sigma_h.
  \end{align*}
\end{remark}

\begin{remark}
  There are many different options for the choice of $\vec{w}_h$.
  In this work, we will simply use the continuous normal velocity defined by the level set function and assume there is no tangential component to the velocity field.
  It is possible to use a backward difference to extract a normal velocity from the two discrete level set functions $\Phi_h(\cdot,t^*)$ and $\Phi_h(\cdot,t^* - \tau)$ by
  \begin{equation*}
    \vec{w}_h =
    - \frac{1}{\tau} \frac{ \Phi_h( \cdot,t^* ) - \Phi_h( \cdot,t^* - \tau ) }{ | \nabla \Phi_h( \cdot,t^* ) | }
    \frac{ \nabla \Phi_h( \cdot,t^* ) }{ | \nabla \Phi_h( \cdot,t^* ) | }.
  \end{equation*}
\end{remark}

\smallskip

\section{Understanding the method in one dimension}

In this section, we want to give a better interpretation
of the method and its limitations, using a simple 1D example.
We consider a point moving along a one-dimensional axis with
positive speed $w$,
i.e.\ $\Phi( x, t ) = x - w t$, $x \in \R, t \in [0,T]$.
Furthermore, we consider a time step with $t^* = \tau = T$.

Let ${U}$ be the interval $[ x_0, x_N ]$ and $\tilde{\T}_h$ be its
decomposition into $N$ sub-intervals $\{ e_j \}_{j=1}^N$ of width
$h=\frac{wT}N$, with $e_j = [ x_{j-1}, x_j ]$.
Let $\Gamma_h^{\old} = x_0$ and $\Gamma_h = x_N$, such that
$D_h = U$ and $\hat{\T}_h = \T_h = \tilde{\T}_h$.
Discretizing using Scheme~1 and a piecewise constant discrete space
$V_h$, i.e.\ polynomial degree $k=0$, we obtain a finite volume type
scheme and equation \cref{eq:scheme} simplifies to
\begin{equation*}
  \bigl[ u_h \varphi_h \bigr](x_N) + \frac{\tau}{\gamma} \bigg( \sum_{i=1}^{N-1}
  u_h \rvert_{e_i} \,
  \bigjump{\varphi_h}(x_i) \, w
  \bigg) = u_h^{\old} \, \varphi_h(x_0)  \quad \text{ for all } \varphi_h \in V_h.
\end{equation*}

By fixing the basis of $V_h$ which consists of characteristic functions
$\varphi_j = \chi_{e_j}$, $j=1,\ldots,N$, we obtain a discrete system.
Denoting the vector of unkowns by $\mathsf{u}$ and supposing
that $u_h^{\old}$ is a given scalar $\mathsf{u}^{\old}$ at $x_0$,
we wish to find $u_h = \sum_{j=1}^N \mathsf{u}_j \varphi_j$ with
\vspace{-0.25\baselineskip}
\begin{align*}
  w\gamma^{-1} \mathsf{u}_1 & = \tau^{-1} \mathsf{u}^{\old} \\
  w\gamma^{-1} \left( \mathsf{u}_j -\mathsf{u}_{j-1} \right) & = 0 && \mbox{ for } j = 2, \ldots, N-1 \\
  \tau^{-1} \mathsf{u}_N + w\gamma^{-1} \left( - \mathsf{u}_{N-1}
  \right) & = 0.
  \\[-0.15\baselineskip]
  \intertext{%
    This system is uniquely solvable, yielding a piecewise constant
    solution $u_h$ with
  }
  \\[-1.35\baselineskip]
  \mathsf{u}_j & = \gamma(\tau w)^{-1} \mathsf{u}^{\old} && \mbox{ for } j = 1, \ldots, N - 1 \\
  u_h \rvert_{\Gamma_h} = \mathsf{u}_N & = \mathsf{u}^{\old}.
\end{align*}
For implementation reasons one might want to compute on a domain
larger than $D_h$, say a computational domain $\Omega \supset
D_h$.
For the domain $D$ in \cref{fig:domain-ex}, for example,
$D_h$ is not easily reconstructed without using knowledge of all
intermediate curves.
As we will illustrate now, the construction of the scheme unfortunately
does not immediately carry over to this case.
Unique solvability requires that $\Omega$ resolves $\Gamma_h \cap D_h$
as a domain boundary.
Furthermore, if $\Omega$ does not meet this requirement, we cannot
guarantee mass conservation any more.

We take a slightly larger domain with the same mesh size and extend
the domain by $3h$ --- $\frac 3 2 h$ to the left and $\frac 3 2 h$ to the
right. This means that $x_{0+\frac 3 2} = \Gamma_h^{\old}$ and
$x_{N-\frac 3 2} = \Gamma_h$, i.e.\ $\Gamma_h^{\old}$ and $\Gamma_h$
lie in the inner part of $e_2$ and $e_{N-1}$, respectively.
The system now changes to
\vspace{-0.25\baselineskip}
\begin{align*}
  w\gamma^{-1} \mathsf{u}_1 & = 0 \\
  w\gamma^{-1} \mathsf{u}_2 & = -\tau^{-1} \mathsf{u}^{\old} \\
  w\gamma^{-1} \left( \mathsf{u}_j - \mathsf{u}_{j-1} \right) & = 0 && \mbox{ for } j = 3, \ldots, N-2 \\
  \tau^{-1} \mathsf{u}_{N-1} + w\gamma^{-1} \left(\mathsf{u}_{N-1} - \mathsf{u}_{N-2} \right) & = 0 \\
  - w\gamma^{-1} \mathsf{u}_{N-1} & = 0,
  \\[-0.15\baselineskip]
  \intertext{%
    which is underdetermined, due to the pure Neumann boundary
    conditions. Analogous to pseudo time stepping, we regularize
    the system with an additional mass term of order $\epsilon > 0$,
    the limit $\epsilon \rightarrow 0$ yields the solution
  }
  \\[-1.35\baselineskip]
  \mathsf{u}_1 & = 0 \\
  \mathsf{u}_j & = \gamma(\tau w)^{-1} \mathsf{u}^{\old} && \mbox{ for } j = 2, \ldots, N - 2 \\
  u_h \rvert_{\Gamma_h} = \mathsf{u}_{N-1} & = \gamma(\tau w + \gamma)^{-1} \mathsf{u}^{\old} \\
  \mathsf{u}_N & \rightarrow \infty.
\end{align*}
For $\mathsf{u}_{N}$, the solution is not well defined and diverges in
the $\epsilon$--limit. For $\mathsf{u}_{N-1}$, we observe that mass
conservation is not fulfilled anymore and for the special case $\gamma =
\tau w$ we obtain $u_h \rvert_{\Gamma_h} = \mathsf{u}_{N-1} = \frac 1 2
\mathsf{u}^{\old}$. Thus it is necessary to resolve $\partial D$
with the computational domain $\Omega$.

\section{Numerical results}

\begin{figure}[tb]
  \centering
  \includegraphics[width=0.25\textwidth]{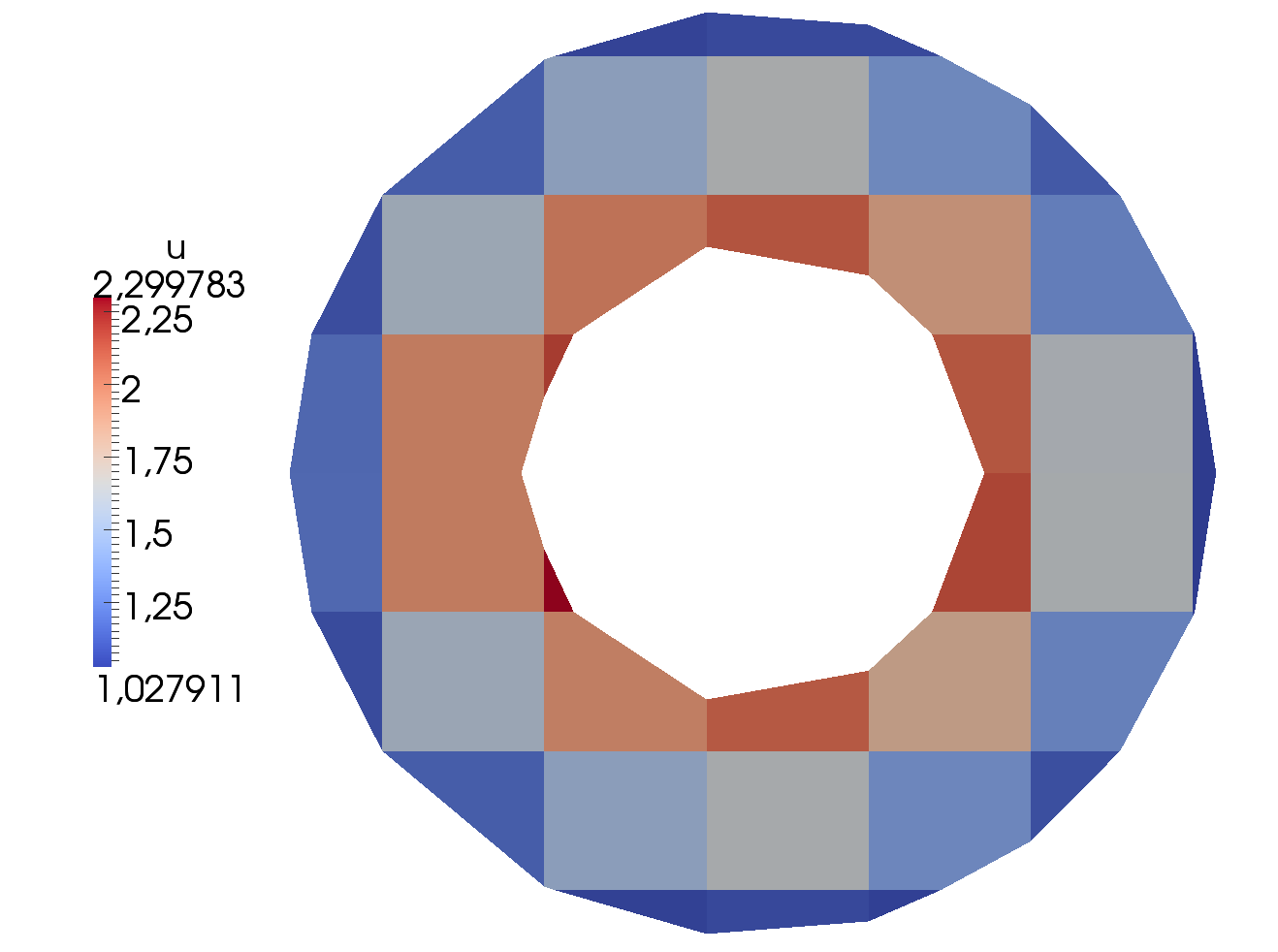}\quad
  \includegraphics[width=0.25\textwidth]{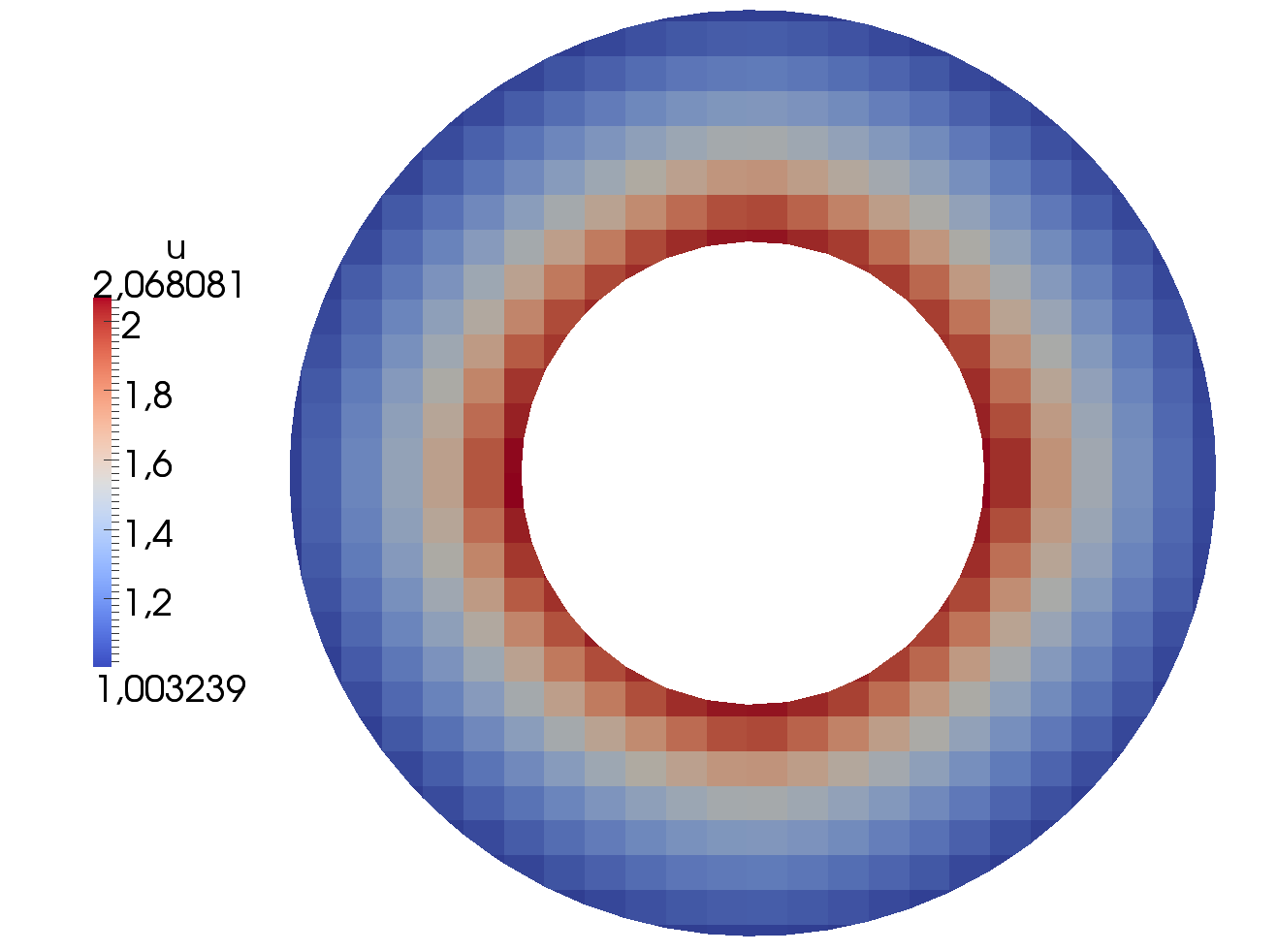}\quad
  \includegraphics[width=0.25\textwidth]{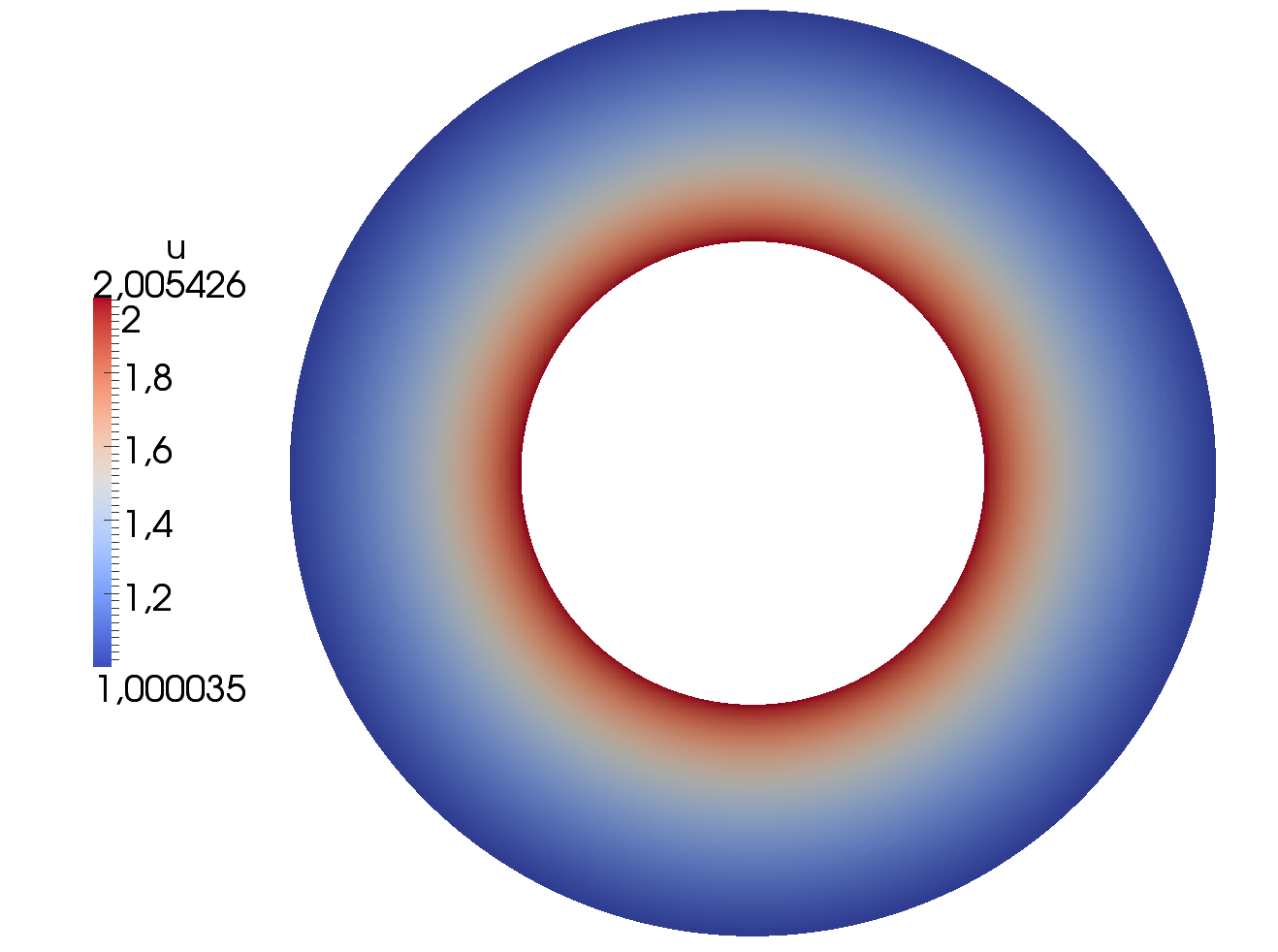}
  \\[2ex]
  \includegraphics[width=0.25\textwidth]{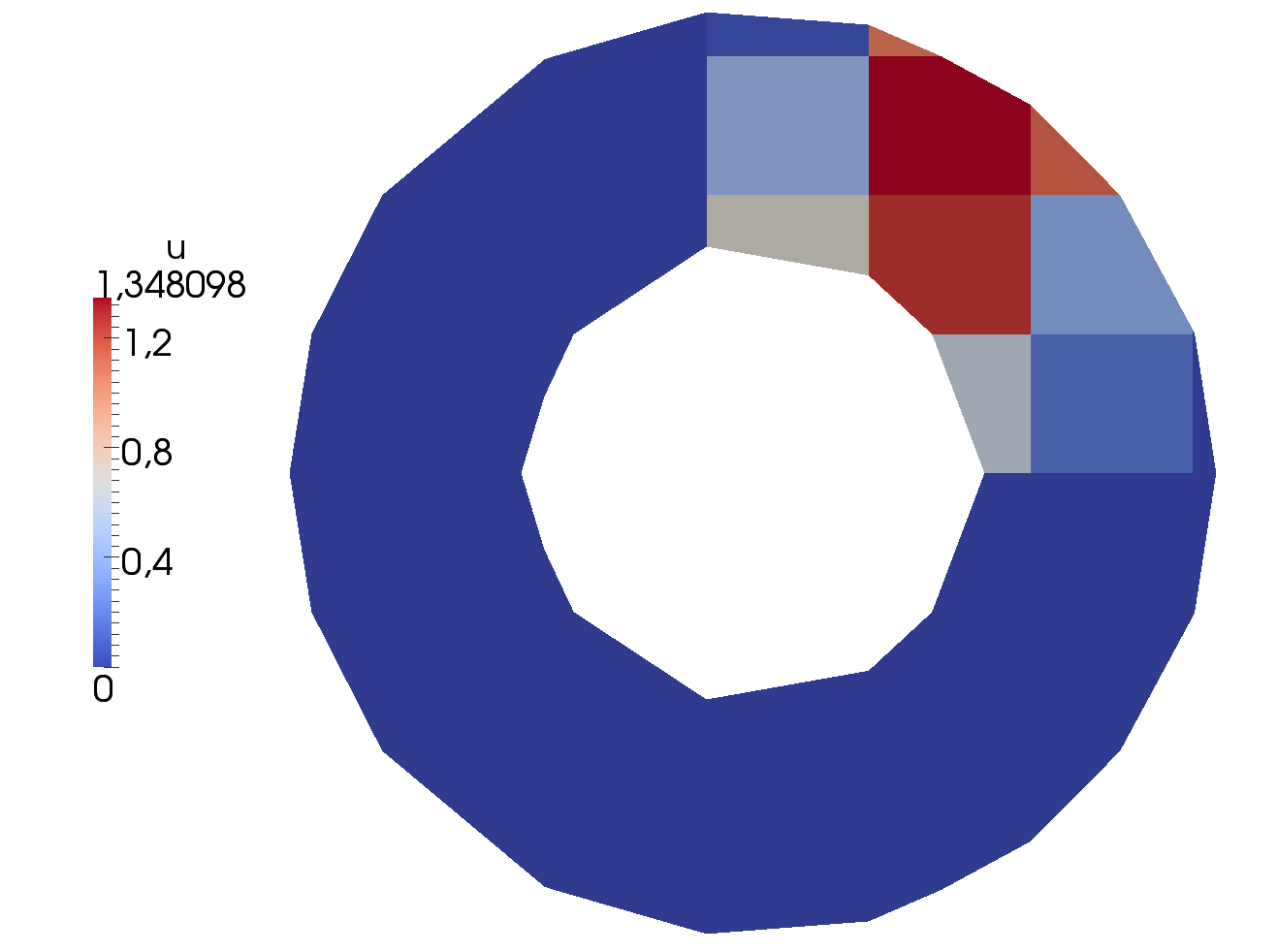}\quad
  \includegraphics[width=0.25\textwidth]{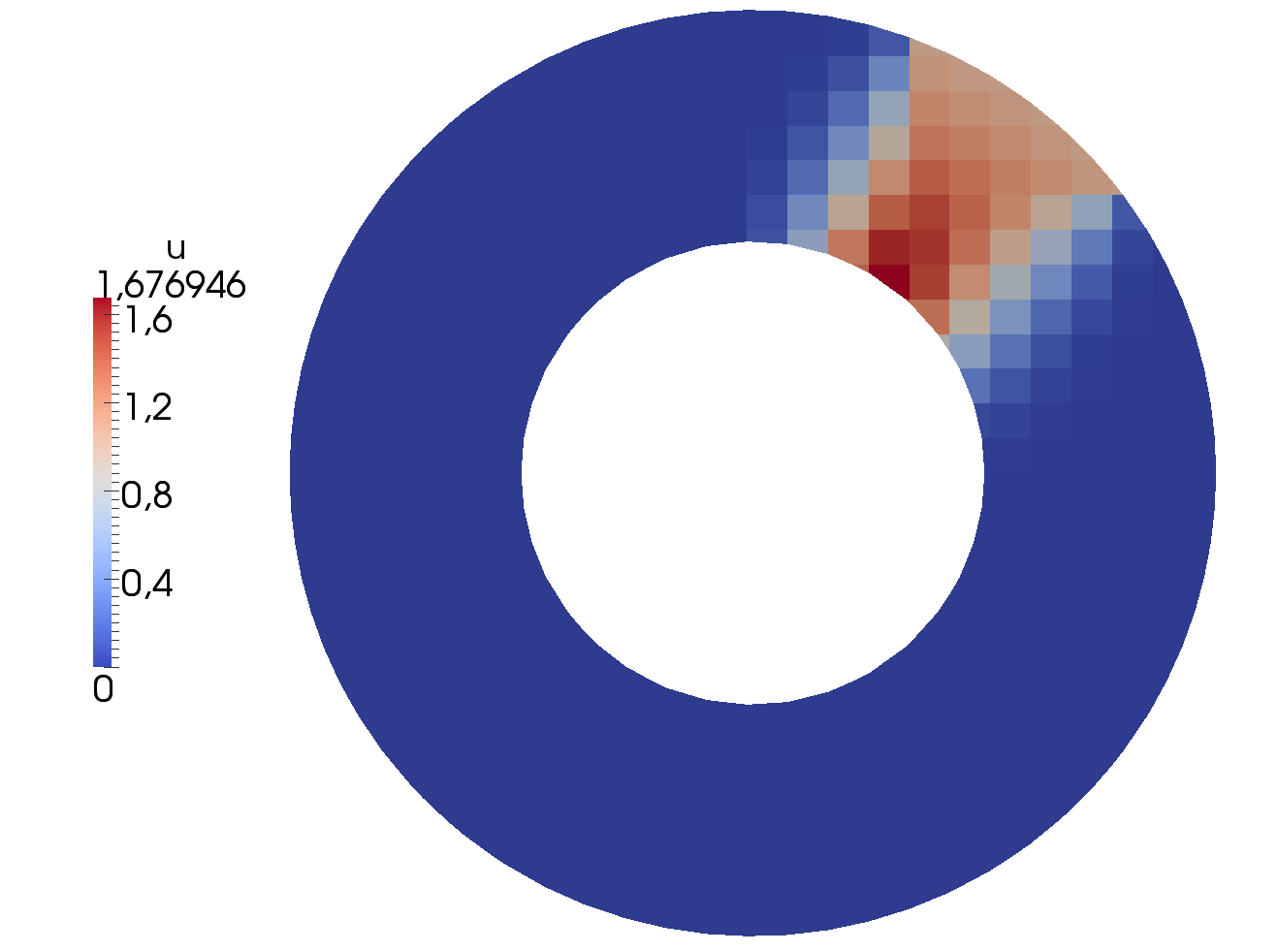}\quad
  \includegraphics[width=0.25\textwidth]{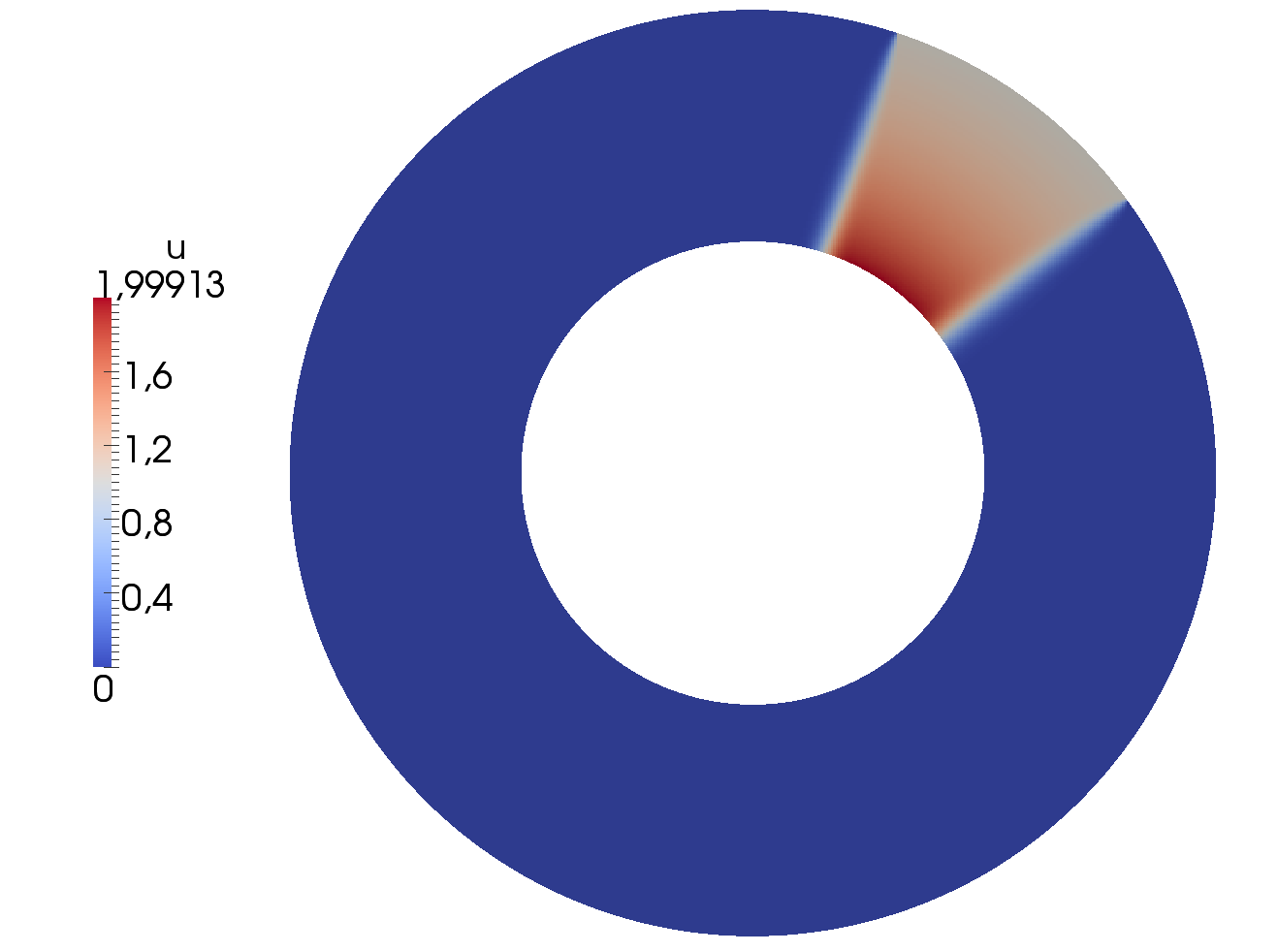}
  \caption{
    Computational domain $D_h$ (discrete level set based reconstruction
    of $D$) and numerical solution $u_h$ under $h$-refinement, using
    $t^* = \tau = \gamma = 0.5$.
    The depicted sequence of numerical solutions in the top row corresponds
    to the cases $\text{\#ncells} = 10^2, 40^2, 640^2$ in
    \cref{tab:annulus_solution_h_refinement}.
    The bottom row considers the same cases for a binary distributed
    $u_h^{\old}$, showing the effect of numerical diffusion.
  }
  \label{fig:annulus_solution_h_refinement}
\end{figure}

The presented scheme has been implemented in the DUNE framework
\cite{dune08:1,dune08:2} using the dune-UDG library \cite{EngHei11}.
Besides providing unfitted DG spaces like the discrete space $V_h$
from \cref{sec:method}, \mbox{dune-UDG} enables the evaluation
of integrals over cut cells $\hat{K} \in \hat{\T}_h$ and their faces
$\hat{E} \in \hat{\E}_h$, which is needed for the assembly of system
matrices in unfitted DG schemes.

The fundamental mesh $\tilde{\T}_h$ (see \cref{sec:discrete_geometry})
used in our code is a structured, Cartesian, quadradrilateral
mesh over a freely choosable domain $U$.
In the following numerical experiments, we use the two-dimensional
domain $U = [-1.5,2] \times [-1.5,1.5]$, decomposed into an equal number
of quadrilaterals in $x$- and $y$-direction.
Furthermore, we employ spaces $V_h$ of polynomial degree $k = 0$.
With respect to an a priori known analytical solution $u$ on
$\Gamma(t^*)$, we compute errors
$\lVert u(t^*) - u_h \rVert_{L^1(\Gamma_h)}$,
$\lVert u(t^*) - u_h \rVert_{L^2(\Gamma_h)}$ and
$\lVert u(t^*) - u_h \rVert_{L^\infty(\Gamma_h)}$.

\subsection{A circle shrinking with constant normal speed}
\label{subsec:annulus_solution_h_refinement}

We consider a circle $\Gamma(t)$ of initial radius $1$, centered at the origin,
which is shrinking with constant normal velocity $\vec{w} \cdot \vec{\nu} = -1$.
Using the description as a level set function, $\Gamma(t)$ can be described
as the zero level set of $\Phi(x,t) = | x | - (1 - t)$, $x \in \R^2$,
$t \in [0,T]$, $T < 1$. Note that $\Phi$ is a signed distance function,
i.e.\ $| \nabla \Phi | \equiv 1$.
Considering a time step with $t^* = \tau = 0.5$, the spatial-only domain $D$
from \cref{sec:motivation} takes the form $D = \{ x \in \R^2 : 0.5 \leq | x | \leq 1 \}$.
Its discrete reconstruction $D_h$ for different values of $h$ can be seen in
\cref{fig:annulus_solution_h_refinement}.

\begin{table}
  \caption{Convergence results under $h$-refinement, using $t^* = \tau = \gamma = 0.5$.}
  \begin{center}
    \footnotesize
    \vspace{-\baselineskip}
    \resizebox{\textwidth}{!}
      { \input{tables/annulus_h_refinement_tau=0.5_results.tex} }\\[-0.82ex]\mbox{}%
  \end{center}
  \label{tab:annulus_solution_h_refinement}
\end{table}

Results for $h$-refinement using $t^* = \tau = \gamma = 0.5$,
$u_h^{\old} \equiv 1$ on $\Gamma_h^{\old}$ and a corresponding
analytical solution $u(t^*) \equiv 2$ on $\Gamma(t^*)$ are
shown in the top row of \cref{fig:annulus_solution_h_refinement}
and in \cref{tab:annulus_solution_h_refinement}.
Note that $\text{\#ncells}$ is the number of cells in the fundamental
mesh $\tilde{\T}_h$.
It does not refer to the number of cut cells $\hat{K} \in \hat{\T}_h$.
Furthermore, since the scheme is globally mass conservative up to machine
precision, the mass written in \cref{tab:annulus_solution_h_refinement}
is both the mass of $u_h^{\old}$ on $\Gamma_h^{\old}$ and the mass of $u_h$
on $\Gamma_h$.
Note that it is an indicator for the geometrical error coming from the
reconstruction of $\Gamma_h^{\old}$, which determines the amount of
mass entering the discrete system.
For the continuous problem and $u^{\old} \equiv 1$, the total mass
in the system equals $2\pi \approx 6.283185$.
We observe convergence of order $1$ in the $L^1$-norm and the $L^2$-norm,
limited by geometrical errors for coarse meshes.
For the $L^{\infty}$-norm, the convergence rate is not that clear but also
seems to be approaching order $1$ for small values of $h$.

\subsection{Non-constant initial concentration}

To illustrate the effect of a non-constant $u_h^{\old}$, we compute the same
setup as in \cref{subsec:annulus_solution_h_refinement} but use a
$u_h^{\old}$ with a binary distribution.
As the worst case scenario, we consider value $1$ for
$0.2\pi < \text{angle} < 0.4\pi$ and $0$ else.
In this scenario, the numerical diffusion is expected to have the largest
impact and for small values of $h$ the jump sharpens, see bottom row of
\cref{fig:annulus_solution_h_refinement}.
The numerical diffusion can be further reduced by using higher-order
methods with flux-limiters.

\section{Conclusion}

We presented a new approach to solve the advection problem
driven by the evolution of an evolving surface.
The method shows attractive properties which justify
further investigation. It is mass
conservative and relatively easy to implement. We avoid constructing
space-time meshes or following characteristics by reformulating the
original problem as a classical transport problem on an unfitted
domain.

In future work, we plan to extend the method to allow
computations on $\Omega \supset D$ and high-order shape
functions. Furthermore, we will apply the method to truly
time-dependent problems and will also consider more general
equations with $\vec{q} \neq 0$.

\input{paper.bbl}

\end{document}


\maketitle

\section{A detailed example}

Here we include some equations and theorem-like environments to show
how these are labeled in a supplement and can be referenced from the
main text.
Consider the following equation:
\begin{equation}
  \label{eq:suppa}
  a^2 + b^2 = c^2.
\end{equation}
You can also reference equations such as \cref{eq:matrices,eq:bb} 
from the main article in this supplement.

\lipsum[100-101]

\begin{theorem}
  An example theorem.
\end{theorem}

\lipsum[102]
 
\begin{lemma}
  An example lemma.
\end{lemma}

\lipsum[103-105]

Here is an example citation: \cite{KoMa14}.

\section[Proof of Thm]{Proof of \cref{thm:bigthm}}
\label{sec:proof}

\lipsum[106-114]

\section{Additional experimental results}
\Cref{tab:foo} shows additional
supporting evidence. 

\begin{table}[htbp]
  \caption{Example table}
  \label{tab:foo}
  \centering
  \begin{tabular}{|c|c|c|} \hline
   Species & \bf Mean & \bf Std.~Dev. \\ \hline
    1 & 3.4 & 1.2 \\
    2 & 5.4 & 0.6 \\ \hline
  \end{tabular}
\end{table}

\bibliographystyle{siamplain}
\bibliography{references}

%% file: tables/annulus_h_refinement_tau=0.5_results.tex
\begin{tabular}{ccccccccc}
  \hline
  \rule{0pt}{0.9em}\#cells & $h$ & $L^1$-error  & (eoc) & $L^2$-error & (eoc) & $L^\infty$-error & (eoc) & mass \\
  \hline
  $5^2$ & $9.220 \cdot 10^{-1}$ & $1.147$ & --- & $9.078 \cdot 10^{-1}$ & --- & $1.002$ & --- & $6.10369$  \\
  $10^2$ & $4.610 \cdot 10^{-1}$ & $3.165 \cdot 10^{-1}$ & $1.86$ & $2.199 \cdot 10^{-1}$ & $2.05$ & $2.998 \cdot 10^{-1}$ & $1.74$ & $6.24228$  \\
  $20^2$ & $2.305 \cdot 10^{-1}$ & $1.951 \cdot 10^{-1}$ & $0.70$ & $1.203 \cdot 10^{-1}$ & $0.87$ & $1.041 \cdot 10^{-1}$ & $1.53$ & $6.27299$  \\
  $40^2$ & $1.152 \cdot 10^{-1}$ & $1.065 \cdot 10^{-1}$ & $0.87$ & $7.175 \cdot 10^{-2}$ & $0.75$ & $8.250 \cdot 10^{-2}$ & $0.34$ & $6.28064$  \\
  $80^2$ & $5.762 \cdot 10^{-2}$ & $5.733 \cdot 10^{-2}$ & $0.89$ & $3.925 \cdot 10^{-2}$ & $0.87$ & $6.628 \cdot 10^{-2}$ & $0.32$ & $6.28255$  \\
  $160^2$ & $2.881 \cdot 10^{-2}$ & $3.075 \cdot 10^{-2}$ & $0.90$ & $2.133 \cdot 10^{-2}$ & $0.88$ & $3.716 \cdot 10^{-2}$ & $0.83$ & $6.28303$  \\
  $320^2$ & $1.441 \cdot 10^{-2}$ & $1.592 \cdot 10^{-2}$ & $0.95$ & $1.081 \cdot 10^{-2}$ & $0.98$ & $2.021 \cdot 10^{-2}$ & $0.88$ & $6.28315$  \\
  $640^2$ & $7.203 \cdot 10^{-3}$ & $8.063 \cdot 10^{-3}$ & $0.98$ & $5.489 \cdot 10^{-3}$ & $0.98$ & $1.005 \cdot 10^{-2}$ & $1.01$ & $6.28318$  \\
  \hline
\end{tabular}

%% file: paper.bbl
\providecommand{\noopsort}{} \renewcommand{\noopsort}[1]{}